\documentclass[11pt,b5paper,english]{article}

\usepackage[left=2.5cm,top=2.5cm,right=2cm,bottom=2cm]{geometry}
\usepackage[T1]{fontenc}
\usepackage{enumerate}
\usepackage{amsmath}
\usepackage{amssymb}
\usepackage{color}
\usepackage{amsthm}
\usepackage{graphicx}
\usepackage[utf8]{inputenc}
\usepackage{lmodern}
\usepackage{float}

\theoremstyle{definition}
\theoremstyle{Notation}
\newtheorem{Notation}{Notation}
\newtheorem{example}{Example}
\newtheorem{lemma}{Lemma} 
\newtheorem{proposition}{Proposition}
\newtheorem{theorem}{Theorem} 

\begin{document}
	
	\title{Note on power hypergraphs with equal domination and matching numbers}
	%\author{María José Chávez de Diego $^\mathrm{a}$, Pablo Montero Moreno$^\mathrm{b}$, \\
	%	María Trinidad Villar-Liñán$^\mathrm{b}$}
	
%	\author{María José Chávez de Diego\footnote{Dpto. de Matem\'atica Aplicada I, Universidad de Sevilla,
%			Avda. Reina Mercedes s/n, Sevilla, Spain, e-mail addresses: mjchavez@us.es }, 		
%Pablo Montero Moreno$^\mathrm{b}$, \\
%		María Trinidad Villar-Liñán\footnote{Dpto. de Geometr\'{\i}a y Topolog\'{\i}a,	Universidad de Sevilla,
%			C/ Tarfia s/n, 41012 Sevilla, Spain, e-mail address: pabmonmor1@alum.us.es,  villar@us.es}$}
	\date{}

\maketitle
\vspace*{-1.5cm}
 \begin{center}
 {\large 	María José Chávez de Diego}\\
  {\it Dpto. de Matem\'atica Aplicada I, Universidad de
 	Sevilla,\\
 	Avda. Reina Mercedes s/n, Sevilla, Spain,\\ 
 	  mjchavez@us.es} \\ 
 
 \vspace*{0.3cm} 
 
 {\large Pablo Montero Moreno, \\
 María Trinidad Villar-Liñán}\\
  	{\it Dpto. de Geometr\'{\i}a y Topolog\'{\i}a,
	Universidad de Sevilla,\\
 	C/ Tarfia s/n, 41012 Sevilla, Spain, 
 	\\  pabmonmor1@alum.us.es,  villar@us.es}\\ 
 
\end{center}

	\begin{abstract}
			We present some examples that refute two recent results in the literature concerning the equality of the domination and matching numbers for power and generalized power hypergraphs. In this note we pinpoint  the flaws in the proofs and suggest how they may be mended.

	\end{abstract}
	
	\section{Preliminaries}

	A {\it (finite)  hypergraph} $  H=(V, E)$ consists of a (finite) set $V$  and a collection $E$ of non-empty subsets of $V$. The elements of $V$ are called {\it vertices} and the elements of $E$ are called {\it hyperedges}, or symply {\it edges} of the hypergraph. A $k-${\it uniform } hypergraph is a hypergraph such that each edge consists of $k$ vertices. A simple graph with no isolated vertices is a 2-uniform hypergraph. Two vertices of $ H$, $u$ and $v$,  are {\it adjacent} if there is an edge $e$ such that $u, v, \in e$. The number of edges containing a vertex $v$ is called the {\it degree} of $v$.
	
	Given a hypergraph  $  H=(V, E),$   $D \subset V$  is   a {\it dominating set} of $ H$ if for every $v\in V - D$ there exists $u\in D$ such  that $u$ and $v$ are adjacent. The minimum cardinality of a dominating set of $ H$, $\gamma(H)$, is  its {\it domination number}. A {\it matching} in $ H$ is a set of disjoint hyperedges. The {\it matching number} of $H$, $\nu(H)$, is the maximum size of a matching in $ H$. A subset $T\subset V$ is a {\it transversal} (or a {\it vertex cover}) of $  H$ if $T$ has nonempty intersection with every hyperedge of $ H$. The {\it transversal number} of $  H$, $\tau(H),$ is the minimum size of a transversal of $H$.
	
 According to \cite{DongSohnLiang2020},  a {\it power hypergraph} $H$ is obtained from a graph $G$ by adding at least one vertex to each edge of $G$; thus, every hyperedge of a power hypergraph contains at least one vertex whose degree is one. The {\it generalized power hypergraph} of a simple graph $G$, denoted $H^{k,s}$, is obtained by blowing up each vertex into a $s$-set and each edge into a $k$-set, where $1\leq s\leq \frac{k}{2}$.
	
Clearly, $H^{k,1}$ is the $k$-uniform power hypergraph. Next we recall a useful notation. 
		
	\begin{Notation}{\rm (\cite{KhanFan2015})} 
		Let $G=(V,E)$ be a simple graph. For any $k\geq 3$ and $1\leq s\leq \frac{k}{2}$, the {\rm generalized power hypergraph} of $G$, denoted by $G^{k,s}$, is defined as the $k$-uniform hypergraph with the vertex set $\{ {\bf v}\,: \, v\in V\}\cup\{{\bf e}\,:\, e\in E\}$, and the edge set $\{{\bf u}\cup{\bf v}\cup{\bf e}\,:\, e=\{u,v\}\in E \}$, where ${\bf v}$ is a  $s$-set containing $v$ and ${\bf e}$ is a $(k-2s)$-set corresponding to $e$.  $G^{k,s}$ is called  the {\rm generalized power hypergraph} obtained from $G$. Particularly, for $s=1$, $G^{k,1}$ is the kth-{\rm power hypergraph} of $G$. 
	\end{Notation}

It is not difficult to check that the domination number is not hereditary for power hypergraphs in general.	However,   the next result holds.

	\begin{proposition}\label{prop:igualdades}{\rm (\cite{DongSohnLiang2020})} 
		Let $H=(V,E)$ be a simple graph, for  $k\geq 3,  \, 1\leq s\leq k/2$  natural numbers, let $H^{k,s}$  denote  the generalized power hypergraph obtained from $H$. We get:
		\begin{enumerate}
			\item $\nu(H^{k,s})=\nu(H)$ and $\tau(H^{k,s})=\tau(H)$.
			\item $\gamma(H^{k, \frac{k}{2}})=\gamma(H)$. 
			\item $\gamma(H^{k, s})=\tau(H^{k, s}), $ if  $1\leq s < k/2.$
		\end{enumerate}
		
	\end{proposition}

\section{Counterexamples to two theorems in Dong et al., 2020.}

A main result from \cite{DongSohnLiang2020} stablishes lower and upper bounds for domination number of power and generalized power hypergraph.

In particular, for any connected generalized power hypergraph $\,H^{k,s}\,$ and $\,k\geq 3$, if $1 \leq s < \frac{k}{2}$, we have $$\nu(H^{k,s})\leq \gamma(H^{k,s}) \leq 2 \nu(H^{k,s})$$ and if $s = \frac{k}{2},$  we get   $\gamma(H^{k,\frac{k}{2}}) \leq \nu(H^{k,\frac{k}{2}})$. (See Theorems 2.1 and 3.1 in \cite{DongSohnLiang2020}).	

	Dong et al. \cite{DongSohnLiang2020} also demonstrate that these bounds are sharp and they give characterizations of extremal hypergraphs for them. Namely the two following theorems are stated and proved.

\begin{theorem}\label{th:2.2}{\rm (Theorem 2.2, \cite{DongSohnLiang2020})} For any connected power hypergraph $H$ of rank $r\geq 3$, $\gamma(H)=\nu(H)$ if and only if $H \in {\mathcal H}_1$. 		
\end{theorem}

\begin{theorem}\label{th:3.4}{\rm(Theorem 3.4, \cite{DongSohnLiang2020})}
		For any connected  generalized power hypergraph $H^{k,s}, $	$\gamma(H^{k,s})= \nu (H^{k,s})$ if and only if $H^{k,s}\in {\mathcal H}^{k,s}_1$ for  $1 \leq s < \frac{k}{2}$ or $H^{k,s}\in {\mathcal H}^{k,k/2}_1$ for $s= \frac{k}{2}$.
\end{theorem}

However,  we present some examples that refute these results concerning the equality of the domination and matching numbers for power and generalized power hypergraphs.

Firstly, let us start by considering  the family ${\mathcal H}^{k,s}_1$, for $k\geq 3$ and  $1 \leq s < \frac{k}{2}$, consisting of generalized power hypergraphs obtained from bipartite connected graphs and ${\mathcal H}_1$ ($={\mathcal H}_1^{k,1}$) denotes the set of power hypergraphs obtained from bipartite connected graphs.

We next study separately the cases $1 \leq s < \frac{k}{2}$ and $s = \frac{k}{2}$.

 \subsection{ The case {\boldmath $1 \leq s < \frac{k}{2}$}}

It is not difficult to refute Theorems \ref{th:2.2} and \ref{th:3.4} by considering  the following counterexamples.

\begin{example}\label{ex:non bipartite}
		Let $H = C_p \vee C_q$ denotes the wedge of two cycles joined by a common vertex, $p, q \geq 3$, being $p$ or $q$ an odd number. The following identities hold.

	\begin{itemize}
			\item $\nu(H) = \left \lfloor \frac{p}{2}  \right \rfloor + \left \lfloor \frac{q}{2}  \right \rfloor;$ 
			
			\item $ \tau(H) = \left \lceil \frac{p}{2} \right \rceil + \left \lceil \frac{q}{2} \right \rceil - 1$
	\end{itemize}	
		From Proposition \ref{prop:igualdades}, we get  $\gamma(H^{k,s})=\tau(H^{k,s})=\tau(H)$ for $1\leq s< k/2.$ Set $n=\nu(H)= \left \lfloor \frac{p}{2}  \right \rfloor + \left \lfloor \frac{q}{2}  \right \rfloor,  $ and let us distinguish two cases.
		
		\begin{enumerate}
			\item If $p$ and $q$ are both odd numbers, then $\left \lceil \frac{p}{2} \right \rceil + \left \lceil \frac{q}{2} \right \rceil - 1=n+1$ and, therefore
			$$\gamma(H^{k,s})=n+1=\nu(H)+1=\nu(H^{k,s})+1.$$
			\item If $p$ is even and $q$ is odd, then $\left \lceil \frac{p}{2} \right \rceil + \left \lceil \frac{q}{2} \right \rceil - 1= n$, and hence 
			$$\gamma(H^{k,s})=n= \nu(H^{k,s}).$$
		\end{enumerate} 
	
	Observe that every graph  of  the form $C_p \vee C_q$, being $p$  even and $q$  odd, is not bipartite while  the  co\-rres\-ponding generalized power hypergraph is extremal for the equality of domination and matching numbers. Therefore, if $s=1$ we found  counterexamples for Theorem \ref{th:2.2} and, if $1<s<\frac{k}{2}$, the counterexamples refute Theorem \ref{th:3.4}.
	
	\end{example}

	Next, we point out what seems to be the main confusion in the characterization given by  Dong et al. in \cite{DongSohnLiang2020}.

	The proofs of both Theorems  2.2 and   3.4 in \cite{DongSohnLiang2020} are quite similar for the case $1\leq s< k/2$, hence we  focus on reasoning only Theorem  2.2.	For the necessary condition it is  affirmed that the equality of domination and matching numbers is raised only by bipartite graphs; however, it is also raised for other graphs, as we have seen in Example \ref{ex:non bipartite}. The misunderstanding is produced because they use   König's Theorem as a characterization of  graphs with equal transversal and matching numbers. But the fact is that the bipartiteness is only a sufficient condition, not a necessary one as we check from \cite{Konig}.
	\begin{theorem} {\rm \bf (König's Theorem, \cite{Konig}) } If $G$ is a bipartite graph, then $\tau(G)=\nu(G).$ 
	\end{theorem}
	
	 Let us notice that the graphs $G$ satisfying $\tau(G)=\nu(G)$ are called {\it König-Egerv\'ary graphs} or said to {\it have the König- Egerv\'ary property}. Trivially, bipartite graphs are examples of such graphs.  König-Egerv\'ary graphs have been extensively studied in the literature; see \cite{Bonomo}  and the references there.   
	Hence,  knowing the family of König-Egerv\'ary graphs  would lead us to  complete the correct characterization of power and generalized power  hypergraphs with equal domination and matching numbers.

	We next sketch a proof.

		Let us denote the family ${\mathcal KE}^{k,s}_1$, for $k\geq 3$ and  $1 \leq s < \frac{k}{2}$, consisting of generalized power hypergraphs obtained from König-Egerv\'ary graphs. For any  $H^{k,s} \in \mathcal KE^{k,s}_1,$ by Proposition \ref{prop:igualdades}, we get $\nu(H^{k,s}) =\nu(H)=\tau(H)=\gamma(H^{k,s})$. 
		Conversely, if $H^{k,s}$ is a generalyzed power hypergraph obtained from a graph $H$ and  such that $\nu(H^{k,s}) =\gamma(H^{k,s})$, then by Proposition \ref{prop:igualdades}, $\tau(H^{k,s}) =\gamma(H^{k,s})$  and $\nu(H)=\tau(H)$ hold and, consequently, we reach $H^{k,s} \in \mathcal KE^{k,s}_1.$ 
		
		Different  characterizations of  König-Egerv\'ary graphs are given in terms of forbidden subgraphs (\cite{Bonomo} and references there in).
		Unfortunately, we have not found a description of the family of König-Egerv\'ary graphs  in an analougous way as the families  ${\mathcal G}_1$ and ${\mathcal G}_{\geq 2}$. 
	
	\subsection{ The case {\boldmath $ s = \frac{k}{2}$}} 
	
	The family ${\mathcal H}^{k,k/2}_1$ contains the generalized power hypergraphs obtained from connected graphs of the family ${\mathcal G}_1\cup{\mathcal G}_{\geq 2}$.
	 The family of graphs ${\mathcal G}_{\geq 2}$ was defined by Randerath and Volkmann in \cite{RanderarthVolkmann1999} and  it was depicted  in \cite{Star-UniformGraphs} as the following nine graphs.
	 
	 \begin{figure}[H] 
		\begin{center}
	\includegraphics[width=13cm]{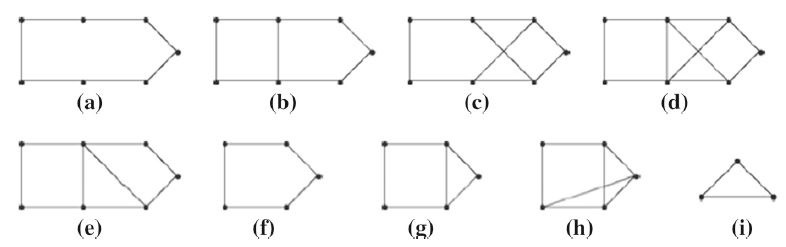}
		\caption{Graphs in $\mathcal{G}_{\geq 2} $ also are included in \cite{DongSohnLiang2020}.}
			\end{center}
		\end{figure}
		
		On the other hand, the family of graphs ${\mathcal G}_1$ was defined by Kano et al. in \cite{Star-UniformGraphs} by using some terminology and notation which are recalled here.
	 The minimum degree of a graph $G$ is denoted by $\delta(G).$   $End(G)$ denotes the set of end-vertices (i.e., vertices of degree one) of $G$. An edge incident with an end-vertex is called a {\it pendant} edge. A vertex adjacent to an end-vertex is called a {\it stem}, and $Stem(G)$ denotes the set of stems of $G$. A graph with a single vertex is called a trivial graph. The {\it corona} $G \circ K_1$ of a graph $G$ is the graph obtained from $G$ by adding a pendant edge to each vertex of $G$. A connected graph $G$ of order at least three is called a generalized corona if $V(G) = End(G) \cup Stem(G)$. \\
		
		Then the family $\mathcal{G}_1$ is that one in which each graph $G$ is the complete graph $K_2$ or a generalized corona, or for each component $G_j$ of $G\setminus(End(G) \cup Stem(G))$ for $j \geq 1$ satisfies one of the following:
		
		\begin{enumerate} [i)]
			\item $G_j$ is a trivial graph.
			\item $G_j$ is a connected bipartite graph with bipartition $V_1$ and $V_2$, where $1 \leq \left| V_1\right| < \left| V_2\right|$. Let $U_{G_j} = V(G_j)\cap N_G(Stem(G))$. Then $\emptyset \neq U_{G_j} \subseteq V_2$ and for any two distinct vertices $x_1, x_2 \in V_1$ that are adjacent to a common vertex of $V_2$, there exist two distinct vertices $y_1,y_2 \in V_2 \setminus U_{G_j}$ such that $N_{G_j}(y_i)=\left\{x_1,x_2 \right\},$ $  i \in \left\{1,2 \right\}$.
			\item $G_j$ is a graph isomorphic to (f), (g), (h) or (i) shown in the precedent figure, and $\gamma(G_j \setminus V_1)=\gamma(G_j) ~ \forall~ \emptyset \neq V_1 \subseteq U_{G_j} \subset V(G_j)$, where $U_{G_j} = V(G_j)\cap N_G(Stem(G))$.
		\end{enumerate}
		 The following example refutes Theorem 3.4 in Dong et al., \cite{DongSohnLiang2020} for $s=\frac{k}{2}.$ 
		\begin{example}	
			Let us consider the  graph $G=K_{2,n}$, $n\geq 2$. Hence $\gamma(G^{k,\frac{k}{2}})=  2 = \nu(G^{k,\frac{k}{2}}).$ However, it is easy to check that $G$ does not belong to $ \mathcal{G}_1\cup \mathcal{G}_{\geq 2} $.     
			
		\end{example}

	An analysis of the proof of Theorem 3.4 in \cite{DongSohnLiang2020} for the case  $s = \frac{k}{2}$   (observe that the case $s = \frac{k}{2}$ does not appear in Theorem 2.2 in \cite{DongSohnLiang2020}) shows that the necessary condition is based on necessary conditions of  Lemmas 3.2 and 3.3 in \cite{DongSohnLiang2020}. 	However, Lemma 3.2 in \cite{DongSohnLiang2020}, as stated, is quite misleading.
	
		We include here the complete characterization of graphs with equal matching and domination numbers   summarized from \cite{DongSohnLiang2020, Star-UniformGraphs, RanderarthVolkmann1999}. 
		
\begin{lemma}{\rm (\cite{DongSohnLiang2020, Star-UniformGraphs})} Let $G$ be a connected graph with $\delta(G)=1$. Then $\gamma(G)=\nu(G)$ if and only if $G\in\mathcal{G}_1.$ 	
\end{lemma}		
		 
	\begin{lemma}{\rm (\cite{Star-UniformGraphs, RanderarthVolkmann1999})} Let $G$ be a connected non bipartite graph with $\delta(G)\geq 2$. Then $\gamma(G)=\nu(G)$ if and only if $G\in\mathcal{G}_{\geq 2}.$
			
	\end{lemma}
			
\begin{lemma}{\rm (\cite{Star-UniformGraphs, RanderarthVolkmann1999})} Let $G$ be a connected  bipartite graph with $\delta(G)\geq 2$ and bipartition of vertex set $X\cup Y$ with $|X|\leq |Y|$. Then $\gamma(G)=\nu(G)$ if and only if $G$ possesses the following properties:
	\begin{enumerate}
		\item $\nu(G)=\gamma(G)=|X|;$ 
		\item for any  two distinct vertices $x_1, \, x_2 \in X$ that are adjacent to a common vertex of $Y,$ there exist two distinct vertices $y_1, \, y_2 \in Y$ such that $y_i$ is adjacent precisely to $\{x_1, x_2\},$ for $i=1,2$.
	\end{enumerate}
\end{lemma}
		
As a consequence, by  using Proposition \ref{prop:igualdades},   the characterization of generalized power hypergraphs $H^{k,\frac{k}{2}}$ with equal domination and matching numbers is  immediately deduced.

\vspace*{0.5cm}

\noindent{\bf Aknowledgements.} Authors are very thankful to Dr. Antonio Quintero Toscano for his helpful coments while preparing this note.

%\vspace*{0.5 cm}

\end{document}